\documentclass[11pt]{article}
\usepackage{amsmath, amsthm, amssymb}    % May not all be necessary
\usepackage{graphicx}			% For including pictures
\usepackage{url}
\newtheorem{theorem}{Theorem}
\newtheorem{corollary}[theorem]{Corollary}
\newtheorem{conjecture}[theorem]{Conjecture}
\newcommand{\msize}[1]{{\left|#1\right|}}
\newcommand{\Pa}[1]{\left(#1\right)}

\begin{document}

% -----------------------------------------------------------------------------

\title{Design Schemes for Fair Dice}

%\author{Tomoko Taniguchi\textsuperscript{1} and Ryuhei Uehara\textsuperscript{1}
\author{Tomoko Taniguchi and Ryuhei Uehara}
%Carlo H. S\'equin\textsuperscript{1} and Another Coauthor\textsuperscript{2}
%\vspace{10pt}\\
% \textsuperscript{1}CS Division, University of California, Berkeley; sequin@cs.berkeley.edu\\
% \textsuperscript{2}Another Famous Institution, Luzern, Switzerland; coauthor@google.com} % end \author
%\textsuperscript{1}
\date{School of Information Science, JAIST, Japan; \{tomoko-t,uehara\}@jaist.ac.jp}
% superscripts are only needed if there is more than one author
% \date{[Draft as of \today]}	% uncomment to use for your own draft purposes
%\date{October 25, 2018}% Suppress any date on submissions
% -----------------------------------------------------------------------------

\maketitle

% Prevent page number 1 from being printed on the first page.
\thispagestyle{empty}

\begin{abstract}
A cube is used as a fair die of 6 faces. However, there are many dice of different shapes on the market.
To make them fair, most of them usually have some symmetric shapes.
We here classify these variants of dice on the market into two groups.
We first consider that a sphere as a model of a fair die with $\infty$ faces.
Based on this model, many symmetric shapes can be modeled as dice obtained by caving spheres.
We also have a familiar fair device; a coin.
That is, a fair coin can be seen as a fair die with 2 faces.
However, a real coin has a thickness, and hence it is, in fact, an unfair die with 3 faces.
From this viewpoint, we propose a way for designing a fair die with $n$ faces for arbitrary $n$.
\end{abstract}

% Bridges papers are usually no more than 8 pages in length.  So
% there's really no need to have numbered sections, unless the
% author really needs to refer to sections by number within the paper's text.  
% So to suppress sequential section numbers, append an asterisk to 
% the \section command, as in:

%%%%%%%%%%%%%%%%%%%%%%%%%%%%%%%%%%%%%%%%%%%
\section*{Introduction}

\begin{figure}[h!tbp]
\centering\vspace{10mm}
\hspace*{-10mm}
\includegraphics[width=4cm,bb=0 0 369 210]{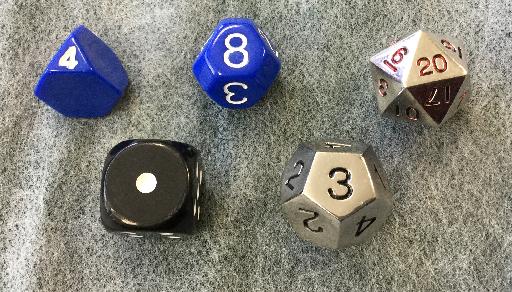}
\hspace{2cm}
\includegraphics[width=4cm,bb=0 0 369 197]{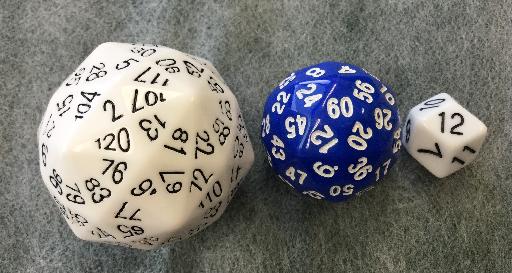}
\caption{Fair dice based on regular and semi-regular polyhedron.}
\label{fig:regular}
\end{figure}

Dice have a long history, and it is uncertain where they originated (see Wikipedia \cite{Wikipedia}).
% It is said that a bone of heel of an animal was used as a die, 
% and this is the first die in history.
Usually, a cube is used as a fair die of 6 faces.
The main reason is that each face of a cube is symmetric, and hence we can conclude that it is \emph{fair}.
The fairness of a die is the most important property, of course.
Intuition tells us that the solid has some symmetric shape, and it can be used as a fair die.
For example, regular polyhedron can be used as a fair die.
Therefore, nowadays, there are many kinds of dice based 
on regular and semi-regular polyhedron on the market (see \figurename~\ref{fig:regular}). 
Unfortunately, however, the possible number of faces is limited when they are 
based on regular and semi-regular polyhedron since the kinds of (semi-)regular polyhedra are 
limited (see Wikipedia \cite{Wikipedia}).
In this paper, we consider two design schemes for designing a fair die for any given number $n$.
We first consider a sphere of uniform density as an ideal fair die.
When you toss a sphere, it will eventually halt at any point with uniform distribution.
That is, a sphere can be considered as  a fair die with $\infty$ faces.
We first discuss a scheme for designing of a fair die with finite faces from a sphere.

We now turn to the other well-known fair device; a fair coin.
In the real world, a coin is frequently used to decide one out of two choices.
That is, when we toss a coin of uniform density, 
we consider that the probabilities of head and tail are equal.
However, a coin has another side out of head and tail in the real world.
Namely, if a coin is thick enough, the probability that side appears 
(or coin stands on its edge) cannot be avoidable.

This thick coin problem had been investigated in the area of probability.
In \cite[Problem 38]{Mos65}, Mosteller asks how thick should a coin be to have a $1/3$ chance of 
landing on edge. He assumed that the coin would slowly settle to its side or its face 
and investigate its outscribed sphere to compute the probability.
However, this mathematical model does not work in the real world.
Sugihara had computed the probability of each face of a standard die according to 
this model\footnote{Sugihara invented his geometric model by himself, 
and independently reached the same model proposed in \cite{Mos65}.},
on the contrary, real experiments with dozen thousands of dice (in a TV show) 
did not fit his evaluation \cite{Sug04a}.

Later, the thick coin problem has been investigated in the area of physics.
In \cite{YM11}, Yong and Mahadevan mentioned that 
the mathematical model (they call ``geometric model'') does not explain 
the real phenomenon when a thick coin is really tossed.
Based on experiments, they propose a ``dynamic model'' which gives us 
a better model of a thick coin to produce probability $1/3$ for each case.
They also analyze their dynamic model and conclude with a certain model.

Now, we turn to the design problem of a fair die of $n$ faces.
We imagine two extreme coins. When it has thickness 0, clearly, 
the coin has no chance of landing on edge. On the other hand, if the coin has a quite huge thickness,
like a colored pencil, it is certain that the coin lands on edge, or it never stands on its face.
Therefore, for any given integer $n$, we should have a reasonable thickness $t(n)$ that gives 
a fair thick coin that stands on its face with probability $2/n$ with $t(3)<t(4)<t(5)<\cdots$.
This argument holds for an $(n-2)$-prism. 
That is, if $(n-2)$-prism has thickness 0, or just a flat regular $(n-2)$-gon, 
it has no chance of landing on edge.
In contrast, if the $(n-2)$-prism is quite long, like a pencil, 
it certainly works as a die of $(n-2)$ faces.
Therefore, for any given integer $n\ge 3$, 
we should have a certain thickness $t'(n)$ of $(n-2)$-prism that is a fair die of $n$ faces,
with $t'(3)<t'(4)<t'(5)<\cdots$.
\begin{figure}[h!tbp]
\centering\vspace{8mm}
\begin{minipage}{0.32\textwidth}
\hspace*{-8mm}
\includegraphics[width=3.5cm,bb=0 0 369 197]{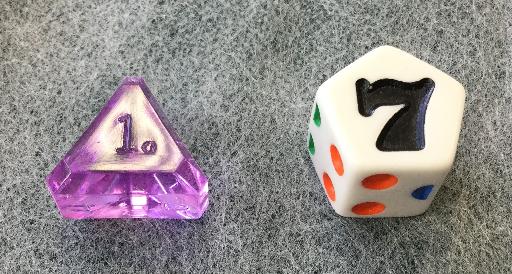}
\caption{Two dice of odd faces on the market.}
\label{fig:5dice}
\end{minipage}
% \end{figure}
\hspace*{8mm}
% \begin{figure}[h!tbp]
\begin{minipage}{0.32\textwidth}
\centering%\vspace{10mm}
%\hspace*{3mm}
\includegraphics[width=3.5cm,bb=0 0 365 183]{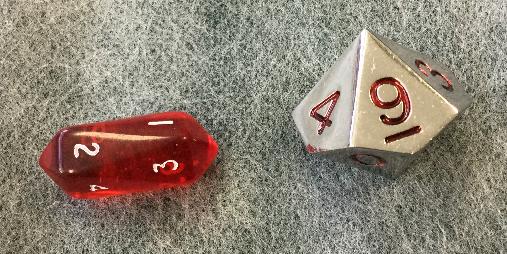}
\caption{Fair, but redundant dice.}
\label{fig:bipyramid}
\end{minipage}
% \begin{figure}[h!tbp]
% \begin{minipage}{0.3\textwidth}
% \centering%\vspace{10mm}
% \hspace*{-10mm}
% \includegraphics[width=3.5cm,bb=0 0 369 216]{hundred.JPG}
% \caption{Two different dice of 100 faces.}
% \label{fig:100dice}
% \end{minipage}
\end{figure}
%
%100$B2sCf!$(B
% 3$B3QCl!'(B51$B2s$NE7(B
% 5$B3QCl!'(B28$B2s$NE7(B
%
%
We first applied the dynamic model in \cite{YM11} for larger $n$
and designed a fair die with $n$ faces. Unfortunately, however,
we found that this model did not fit experimental results for large $n>4$.
In fact, we found two dice of faces 5 and 7 based on the same idea (\figurename~\ref{fig:5dice}), 
nevertheless, one of them is not fair. 
We tossed 100 times and obtained head/tail 51 times for the die of face 5,
which is bigger than the expected value 40.
On the other hand, we obtained head/tail 28 times out of 100 for the die of face 7,
which fits the expected value $200/7=28.6$.
We have not yet succeeded to analyze the reason. 
In fact, the reason why the die of face 5 is not fair may be the rounded corners.
As mentioned by Sugihara \cite{Sug04b},
the center of gravity has a tendency to get out of position by rounding corners.

From the experiments, we observe that when $n$ is bigger, the think coin becomes longer like a log, 
and it rolls when it is tossed. This rolling log phenomenon seems to be the main reason 
why the probability that the log stands on its face is much smaller than the one computed by the dynamic model.
We will show a case study for $n=13$ to design a fair die with $13$ faces based on experiments.

It is worth mentioning that, from a practical point of view, 
when we admit ``redundancy'' to a die, there are two simple ways for designing 
a fair die as follows (\figurename~\ref{fig:bipyramid}):

% \begin{figure}[h!tbp]
% \centering\vspace{10mm}
% \includegraphics[width=4cm,bb=0 0 365 183]{bipyramid.JPG}
% \caption{Fair, but redundant dice.}
% \label{fig:bipyramid}
% \end{figure}

\paragraph{Double Face Method:}
When we allow that a die has two or more faces of the same pip, there is a simple solution. 
For any integer $n>2$, $n$-gonal bipyramid can be used as a fair die
when each number $i$ with $1\le i\le n$ appears twice on faces.
In this paper, we consider this solution is not satisfactory since half faces are redundant.

\paragraph{Unused Face Method:}
As a hexagonal pencil can be used as a die with 6 faces in elementary schools,
the $n$-prism can be used a fair die. 
However, from a theoretical viewpoint, it has a positive probability that 
the pencil stands up on its end.
To avoid this, we can sharpen both sides of the pencil.
In this time, the pencil cannot stand up on the end since it is now just a point of area 0. 
However, in this case, around the point, 
we have some surface area of the pencil that will never be used as a face of the die.
In this paper, we consider such a surface area redundant.

\section*{Sphere Model and its Variants}

In this section, we propose a fair die model based on a sphere of uniform density as an ideal fair die.
We consider a sphere is a fair die with $\infty$ faces.
Based on this model, we can consider a die obtained by the following procedure is fair for a given $n>1$;
\begin{itemize}
\item Put $n$ points $p_1,p_2,\ldots,p_n$ on a surface of sphere $P$ ``uniformly'', 
      which is defined by the following conditions in this procedure.
\item For each $p_i$, we plane the surface of $P$ around $p_i$ to obtain flat face $f_i$ that contains $p_i$ 
      such that
      (1) the center of gravity of $P$ is not changed, and
      (2) $\msize{f_i}=\msize{f_j}$ holds for each $i\neq j$, where $\msize{f_i}$ is the area of $f_i$ on $P$.
\end{itemize}
When we plane until the original convex surface of $P$ is removed, 
the dice based on regular polyhedron in \figurename~\ref{fig:regular} can be explained by this model.
Each die based on semi-regular polyhedron also can be explained by this model by attaching 
some pyramids glued on a surface of a regular polyhedron properly and symmetry.

From this viewpoint, we can consider that the design of a fair die of $n$ faces 
is an application of the Voronoi diagram (the details of Voronoi diagram are omitted here) on a sphere.
That is, we arrange $n$ Voronoi points $p_1,p_2,\ldots,p_n$ on sphere uniformly in the sense that 
(1) each Voronoi region has the same area, and 
(2) each Voronoi point is placed on the geometric center of each Voronoi region.
If we can find this uniform distribution of $n$ points on a sphere, by planing a bit 
(to not change the center of gravity of $P$), we will obtain a fair die of $n$ faces.
This problem had been investigated a long time in geography and geometry,
and some computer programs are available (see, e.g., \cite{BB99,Yam01}).

% \begin{figure}[h!tbp]
% \centering
% \vspace{10mm}
% \includegraphics[width=5cm,bb=0 0 369 216]{hundred.JPG}
% \caption{Two different dice of 100 faces.}
% \label{fig:100dice}
% \end{figure}

\begin{figure}[h!tbp]
\centering\vspace{10mm}
\begin{minipage}{0.45\textwidth}
%\hspace*{5mm}
\includegraphics[width=4.5cm,bb=0 0 369 216]{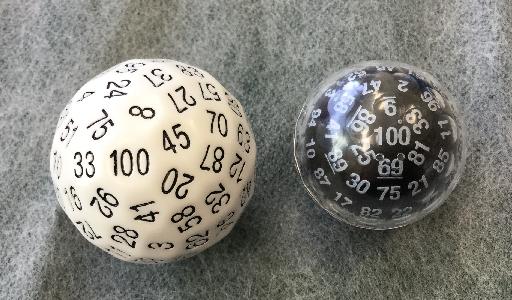}
\caption{Two different dice of 100 faces.}
\label{fig:100dice}
\end{minipage}
\hspace{1cm}
% \end{figure}
% \begin{figure}[h!tbp]
\begin{minipage}{0.45\textwidth}
%\vspace{10mm}
\hspace*{2mm}
\includegraphics[width=4cm,bb=0 0 575 374]{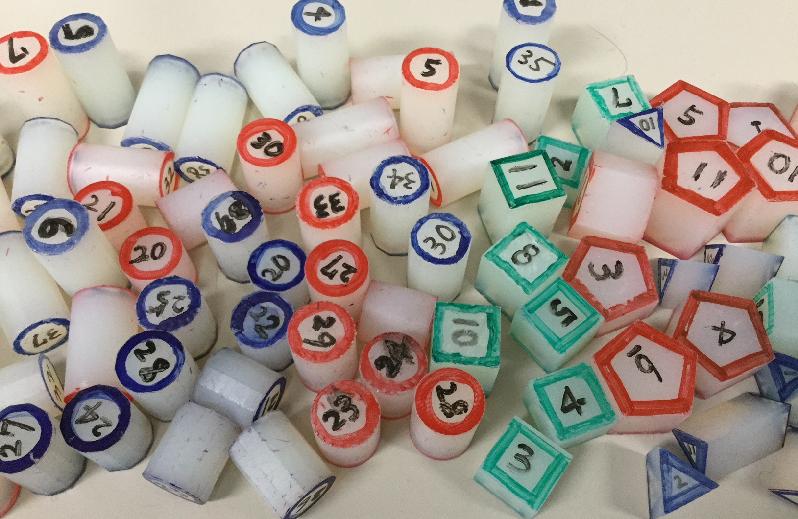}
\caption{Many dice were made by a 3D printer.}
\label{fig:3D}
\end{minipage}
% \begin{minipage}{0.4\textwidth}
% \center
% \includegraphics[width=4cm,bb=0 0 212 228]{golf.JPG}
% \caption{Golf ball as an application of a fair die.}
% \label{fig:golf}
% \end{minipage}
\end{figure}

On the dice with 100 faces on the market, there are two different types (\figurename~\ref{fig:100dice}).
The left one in \figurename~\ref{fig:100dice} has two special numbers 1 and 100, 
which correspond to two special flat faces that are larger than the other.
Therefore, we can conclude that this die should not be fair (as mentioned in Wikipedia \cite{Wikipedia}).
On the other hand, the right one in \figurename~\ref{fig:100dice} can be explained by our model;
each $p_i$ has a dimple of the same size, and the arrangement of these dimples is symmetric, 
which guarantee that the center of gravity of $P$ is not changed\footnote{This die has 
an empty space inside and partially filled with sand. 
By this trick, the die does not roll so long time, but we need to throw carefully to make it fair.}.
Therefore, we can conclude that this die is fair.
This fair die reminds us a golf ball.
We mention that the designing of a fair die with many faces has an application to 
the designing of a good golf ball.

\noindent
\textbf{Note:} After preparing this draft, we found that 
some dice designed based on the same idea in this section are on the market.
See \url{http://mathartfun.com/d357.html} for their products.
%
% Radar dome

\section*{Coin Model and its Extension}
Hereafter, we consider cylinders and $n$-prisms.
We standardize that each cylinder has radius 1, and 
each $n$-prism is inscribed in a circle of radius 1.

In \cite{YM11}, Yong and Mahadevan proposed a dynamic model, which is based on a physical model.
Their model gives us the following theorem:
\begin{theorem}[\cite{YM11}]
\label{th:dynamic}
Let $t$ be the thickness of a coin (of radius 1).
Then the coin stands on edge with 
the probability 
\[
P_D(t)=1-\frac{2}{\pi}\cos^{-1}\Pa{\frac{t}{\sqrt{1+t^2}}}.
\]
\end{theorem}
We aim at designing of a fair die with $n$ faces. To achieve this goal, we consider the case
that the probability $P_D(t)=\frac{n-2}{n}$. Solving this equation, we have the following corollary.
\begin{corollary}
\label{co:dynamic}
In Theorem \ref{th:dynamic}, we let $P_D(t)=\frac{n-2}{n}$ and solve this equation for $t(n)$.
Then we have the following:
\[
t(n)=\msize{\frac{\cos \pi/n}{\sin \pi/n}}.
\]
\end{corollary}
For $n=3,4,\ldots,13$, we can obtain Table~\ref{tab:3-13}.

\begin{table}
\centering
\caption{Thickness of a coin that is on a face with probability $2/n$ based on the dynamic model.}
\label{tab:3-13}
\begin{tabular}{|r|ccccccccccc|}\hline
$n$    & 3 & 4 & 5 & 6 & 7 & 8 & 9 & 10 & 11 & 12 & 13 \\\hline
$t(n)$ & 0.577 & 1.000 & 1.376 & 1.732 & 2.077 & 2.414 & 2.747 & 3.078 & 3.406 & 3.732 & 4.057 \\\hline
\end{tabular}
\end{table}

% \begin{figure}[h!tbp]
% \centering
% \vspace{10mm}
% \includegraphics[width=5cm,bb=0 0 575 374]{dice.jpg}
% \caption{Many dice made by a 3D printer.}
% \label{fig:3D}
% \end{figure}

\subsection*{Cylinder as a Fair Die of 13 Faces}
To check the (in)validness of the dynamic model,
as an example, we consider the extreme case $n=13$,
which is assumed to be used as a kind of a fair die with 13 faces.
First, we consider a cylinder that stands on two circular faces with probability $2/13=0.154$,
and that lays on its side with probability $11/13$.
In Table~\ref{tab:3-13}, the dynamic model tells us to make a cylinder of height $4.057$ for radius $1$.
Therefore, we made many cylinders of radius 10mm of various heights and tossed them 100 times for each (\figurename~\ref{fig:3D}).
The experimental results are summarized in Table~\ref{tab:cylinder}:

\begin{table}
\centering
\caption{The number of heads/tails out of 100 tosses of each cylinder; 
it was expected to be $100\times 2/13=15.4$ when height is $41$mm.}
\label{tab:cylinder}
\begin{tabular}{|r|ccccccccccccccc|}\hline
Heights (mm) & 20 & 21 & 22 & 23 & 24 & 25 & 26 & 27 & 28 & 29 & 30 & 31 & 32 & 33 & 34 \\\hline
Heads/tails   &  9 & 19 &  8 &  5 &  8 &  7 &  5 &  9 &  3 &  3 &  3 &  4 &  3 &  1 &  4 \\\hline\hline
Heights (mm) & 35 & 36 & 37 & 38 & 39 & 40 & 41 & 42 &    &    &    &    &    &    &    \\\hline
Heads/tails   &  3 &  2 &  6 &  1 &  7 &  2 &  2 &  4 &    &    &    &    &    &    &    \\\hline
\end{tabular}
\end{table}

%$B:G=i$N<B83$N0J2<$OL5:n0Y$K<N$F$?(B
%35 36 37 38 39 40
%3  3  3  1  4  0

The experiments above were done with the following condition: 
they were tossed from a dish from height 80cm on a carpet.
We also made experiments of different heights (20cm/40cm/80cm) and 
different floors (carpet/wooden desk), but the behavior is not changed.

\subsection*{Design of a Fair Die of 13 Faces based on Cylinder Model}
For any integer $n\ge 3$, let $P_n$ be a regular $n$-gon.
We also let $C_n$ and $c_n$ be the outscribed circle and the inscribed circle to $P_n$, respectively.
Without loss of generality, we suppose that $C_n$ is a circle of radius 1.
That is, all $C_n$s are congruent with a circle $C$ of radius 1.
Then the radius $r_n$ of $c_n$ can be computed by $\cos\frac{\pi}{n}$.
For given $P_n$ and height $h$, 
we denote by $Pr_{HT}(P_n,h)$ the probability that 
the $n$-prism of height $h$ stands on its two faces of regular $n$-gon when it is tossed.
For the circles $C$ and $c_n$, the probabilities $Pr_{HT}(C,h)$ and $Pr_{HT}(c_n,h)$ are defined in the same manner.
Then we conjecture the following:
\begin{conjecture}
\label{conj:in-out}
\[
Pr_{HT}(c_n,h)<Pr_{HT}(P_n,h)<Pr_{HT}(C,h).
\]
\end{conjecture}
For $n=3,4,5,11$, we made several experiments, and Conjecture~\ref{conj:in-out} follows that.
We also confirmed that $n>4$, the dynamic model does not fit to experiments.

Based on Conjecture \ref{conj:in-out}, we designed a fair die of 13 faces by a regular $11$-prism.
In this case, $r_{11}=0.9595$ for the inscribed circle to $P_{11}$.
We first tossed several $11$-prisms and cylinders 500 times as in Table~\ref{tab:fair11}.

% $B1_Cl(B
% $B9b$5(Bmm 20 21 22 23 24 25 26 27 28 29
% $BE7CO(B400$B!!(B   73 60 45 59 42 36 42 39 29 29
% $BE7CO(B100$B!!(B   19 13 22 13  8 16 11 10  9  9
% $BE7CO(B500$B!!(B   92 73 67 72 50 52 53 49 38 38

% 11$B3Q7A!!(B
% $B9b$5(Bmm 20 21 22 23 24 25 26 27 28 29
% $BE7CO(B400   68 54 52 46 29 22 33 16 22 21
% $BE7CO(B100$B!!(B 16 20 24  8  7  8  6  3  8  4
% $BE7CO(B500$B!!(B 84 74 76 54 36 30 39 19 30 25

\begin{table}
\centering
\caption{The number of heads/tails out of 500 tosses of each cylinder and each 11-gon.
They were tossed from height 20cm to on a wooden desk.
When they are estimated as a fair die of 13 faces, it is expected to produce 
$500\times 2/13=76.92$ times.}
\label{tab:fair11}
\begin{tabular}{|r|cccccccccc|}\hline
Heights (mm)              & 20 & 21 & 22 & 23 & 24 & 25 & 26 & 27 & 28 & 29 \\\hline
Heads/tails of cylinders & 92 & 73 & 67 & 72 & 50 & 52 & 53 & 49 & 38 & 38 \\
Heads/tails of $n$-gons  & 84 & 74 & 76 & 54 & 36 & 30 & 39 & 19 & 30 & 25 \\\hline
\end{tabular}
\end{table}

By fitting a linear function $t(h)=ah+b$, where the number of heads/tails $t(h)$ out of 500 tosses 
for height $h$, we obtain $a=-5.30909$ and $b=188.473$ for the cylinders
and $a=-7.12121$ and $b=221.17$ for the 11-prisms.
Then we can observe that Conjecture \ref{conj:in-out} holds for $n=11$:
When we let the probability to be $2/13$, 
the corresponding heights of the cylinder of radius $r_{11}=0.9595$, the 11-prism, 
and the cylinder of radius 1 are $20.160$, $20.255$, and $21.011$, respectively.
These results support Conjecture \ref{conj:in-out}, and we conclude 
that the fair die of 13 faces based on the 11-prism is of height 20.3mm, which is between 20.2mm and 21.0mm.

% f(x)=a*x+b
% fit f(x) "7.dat" using 1:2 via a,b
% a               = -5.30909         +/- 0.7146       (13.46%)
% b               = 188.473          +/- 17.63        (9.353%)

% fit f(x) "7.dat" using 1:3 via a,b
% a               = -7.12121         +/- 1.124        (15.79%)
% b               = 221.17           +/- 27.74        (12.54%)

% $BE7$,=P$?2s?t(B
% $B!!!!!!!!(B 1$B!!(B 2$B!!(B 3$B!!(B 4$B!!(B 5$B!!(B 6$B!!(B 7$B!!(B 8$B!!(B 9$B!!(B 10$B!!(B 11
% 3$B3QCl!'(B 63  43  35  22  18  10   6   4   2    6    4
% 4$B3QCl!'(B 82  81  71  59  40  31  30  12  16    6    6
% 5$B3QCl!'(B 81  84  66  67  63  51  48  36  32   24   19

% A> $B>e$N?t;z$,Dc$$=g$N9b$5(B
% B> $BK\Ev$N9b$5(B(mm)
% C> $B2<$N?t;z$,(B100$BJ,$NE7CO$,=P$?2s?t(B

% $B$H$7$F!$(B

% 11$B3Q7A!!(B
% A>    1    2    3    4    5    6
% B> 38.5 39.0 39.5 40.0 40.5 41.0
% C>    2    2    4    2    2    0

% 100$B2s?6$j$d$7$?!#=87W$7$d$7$?!*(B

% 11$B3Q7A!!(B
% $B9b$5(Bmm 30 31 32 33 34 35 36 37 38 39 40
% $BE7CO!!!!(B3  8  4  2  2  5  1  2  1  2  3

% $B$I$&$+$7$i$s!#(B
% $B6=L#?<$$7k2L$K$J$C$?$+$7$i$s!#(B

% $B%5%$%3%m(B100$B2s?6$j$^$7$?!#(B

% 11$B3Q7A(B
% $B9b$5(Bmm  20 21 22 23 24 25 26 27 28 29
% $BE7CO!!(B   $B!!(B8 10 10  4  12   6  5 10   7  4

% tomoko6,7$B$O$^$@F~$l$F$J$$(B

\section*{Conclusions}

%rolling log

\begin{figure}[h!tbp]
\centering\vspace{10mm}
\begin{minipage}{0.4\textwidth}
\includegraphics[width=4cm,bb=0 0 369 180]{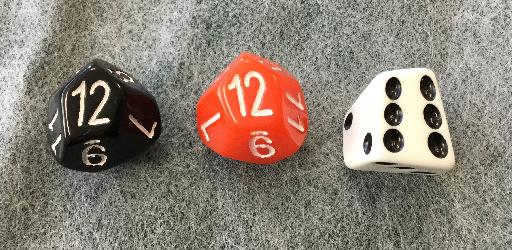}
\caption{Deformed, but fair dice.}
\label{fig:deform}
\end{minipage}
\hspace{10mm}
\begin{minipage}{0.4\textwidth}
\includegraphics[width=3.2cm,bb=0 0 343 227]{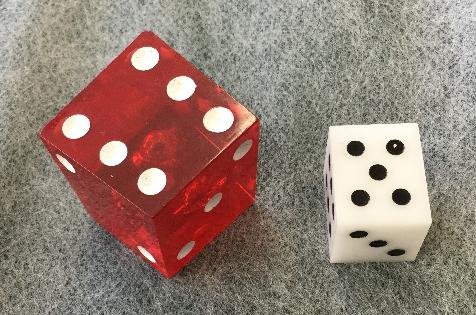}
\caption{Deformed and unfair dice.}
\label{fig:deform-unfair}
\end{minipage}
\end{figure}

In this paper, we classify and propose two models for designing a fair die.
However, there are some dice on the market which cannot be handled in this framework.
The dice in \figurename~\ref{fig:deform} are deformed, but still fair.
Their faces are symmetric, and they are market products as fair dice.
We do not yet classify how can we decide they are symmetric enough to be fair.

% \begin{figure}[h!tbp]
% \centering\vspace{10mm}
% \includegraphics[width=4cm,bb=0 0 343 227]{deform-unfair.JPG}
% \caption{Deformed and unfair dice.}
% \label{fig:deform-unfair}
% \end{figure}
The other dice in \figurename~\ref{fig:deform-unfair} are deformed and unfair.
The left die is totally deformed, and the right die is rectangular.
Especially, analyze and describe the probability of each face of 
the rectangular die of size $a\times b\times c$ seems to be an interesting topic.

The other topic that we did not investigate in this paper is that 
the arrangement of $n$ numbers on a fair die of $n$ faces.
The ordinary cubic die satisfies that summation of two opposite faces makes 7,
and hence there are two different cubic dice such that they are mirror symmetry with each other.
General die of $n$ faces for large $n$, there are many ways to arrange these pairs $(i,n-i+1)$.
The balanced arrangement problem is investigated in \cite{BFS16}.
Making a fair die by caving or/and drilling many holes on the surface, 
such balanced arrangements are another problem to be considered.

\section*{Acknowledgements}

The authors thank Prof. Kokichi Sugihara for his helpful discussions and suggestions.

\raggedright				% no right justification for References

\end{document}